\documentclass{amsart}
\usepackage{amssymb,amsthm, amsmath, amsfonts}
\usepackage{graphics}
\usepackage{hyperref}
\usepackage[all]{xy}
\usepackage{enumerate}
\usepackage[mathscr]{eucal}
\usepackage{bbm}
\usepackage{mathrsfs,amssymb}

\theoremstyle{plain}
\newtheorem{theorem}{Theorem}[section]
\newtheorem{lemma}[theorem]{Lemma}
\newtheorem{proposition}[theorem]{Proposition}

\numberwithin{equation}{section}

\theoremstyle{definition}

\newtheorem{definition}[theorem]{Definition}

\begin{document}
\title[Representation dimension of smash products]{On the representation dimension of smash products}
\author[Zheng ET AL.]{Lijing Zheng$^{1,a}$}
\author[]{Chonghui Huang$^{1,b}$}
\author[]{Qianhong Wan $^{2}$ }


\thanks{$^A$. This work is partly supported by Natural Science Foundation of China \#11271119.\\
$^B$. This work is partly supported by Natural Science Foundation of China \#11201220.\\
E-mail: zhenglijing817@163.com, huangchonghui@usc.edu.cn,  77927023@qq.com. }
\dedicatory{$^1$ School of Mathematics and Physics \\ University of South China \\ Hengyang, 421001, Hunan, P. R. China\\
   $^2$ School of Mathematics and Statistics\\ Hunan University of Commerce\\ Changsha, 410205, Hunan, P. R. China  }

\begin{abstract}
Let $A$ be a finite dimensional $G$-graded algebra with $G$ a finite group, and $A\# k[G]^{\ast}$ be the smash product of $A$ with the group $G$. Our results can be stated as follows: (1) If $A$ is a self-injective algebra and separably graded, then the dimensions of triangulated categories $\underline{\rm mod}A$ and $\underline{\rm mod}A\# k[G]^{\ast}$ are equal. In particular, we obtain that the representation dimension of $A\# k[G]^{\ast}$ is at least the dimension of triangulated category  $\underline{\rm mod}A$ plus 2; (2) Generally, if $A$ is a $k$-algebra and separably graded, then the Oppermann dimensions of  $A$ and $A\# k[G]^{\ast}$ are equal. In particular, we obtain that the representation dimension of $A\# k[G]^{\ast}$ is at least the Oppermann dimension of $A$ plus 2. In the end, we give two examples to illustrate our results.
\end{abstract}

\maketitle

\tableofcontents

\section{Introduction}
\subsection{Background}
The representation dimension of a finite dimensional algebra  was introduced by Auslander \cite{Aus}. Auslander has shown that an algebra is representation finite if and only if its representation dimension is at most two. He expected that the representation dimension should measure how far an algebra is far from being of representation finite algebras.
For more than three decades, it was unclear whether Auslander's philosophy works.

 In 2006, Rouquier showed that the representation dimension of the exterior algebra of an $n$-dimensional vector space is $n+1$, in which the notion of the dimension of a triangulated category played a very important role (cf. \cite{Ro2}).
 In \cite{EHSST}, Erdmann etc. introduced finite generation assumption ($\mathbf{Fg}$). Bergh showed that if a non-semisimple self-injective algebra satisfies $\mathbf{Fg}$, then the dimension of its stable module category as a triangulated category is at least its complexity minus one \cite{Ber}.

In 2009, Oppermann gave a lattice criterion of computing the lower bound of the representation dimesion of a finite dimensional algebra and in this way one obtained a large class of algebras with large representation dimensions \cite{Opp}. For a finite dimensional algebra $A$, the key tool he invented for computing the representation dimension is nowadays called Oppermann dimension (see the definition in Section 2), denoted it by Odim$A$. And he obtained that rep.dim$A$ $\geq$  Odim$A$+2 \cite{Opp}.

On the other hand, representations and homological properties of smash product algebras (see the definition in Section 3) have been widely studied, see \cite{JJ,WZ}. It has been shown that $A\#k[G]^*$ and $A$ share many common properties. For example, $A\#k[G]^*$ is a self-injective algebra if and only if so is $A$ \cite{JJ}. It is natural to compare the dimensions of the stable module categories of $A\#k[G]^*$ and $A$ as triangulated categories.  Concerning the representation dimension of a finite dimensional algebra, it is reasonable to compare the Oppermann dimensions of $A\#k[G]^*$ and $A$. Then it is possible to obtain a new class of algebras whose representation dimensions are arbitrarily large.

Our main results in this paper are the following two theorems.

\begin{theorem}\label{main-theorem} Let $A=\oplus_{g\in G}A_{g}$ be a finite dimensional $G$-graded self-injective $k$-algebra with $G$ a finite group. Then
$${\rm dim}(\underline{\rm mod}A\#k[G]^*)\leq {\rm dim}(\underline{\rm mod}A).$$
Furthermore, if $A$ is separably graded, then
$${\rm dim}(\underline{\rm mod}A\#k[G]^*)={\rm dim}(\underline{\rm mod}A).$$
In this case, we obtain that ${\rm rep.dim}(A\#k[G]^*)\geq {\rm dim}(\underline{\rm mod}A)+2. $

\end{theorem}
\begin{theorem} \label{main-theorem-1}Let $A=\oplus_{g\in G}A_{g}$ be a $G$-graded  $k$-algebra with $G$ a finite group .
Denote by $B=A\#k[G]^*$.  Then $${\rm Odim}B\leq {\rm Odim}A.$$
Furthermore, if A is separably graded, then $${\rm Odim}B={\rm Odim}A.$$
In this case, we obtain that ${\rm rep.dim}B\geq {\rm Odim}A+2.$

\end{theorem}

\subsection{Organization}
This paper is organized as follows. In Section 2, we shall recall some basic definitions and facts needed in our proofs. Theorem \ref{main-theorem} and \ref{main-theorem-1} will be proved in Section 3. In Section 4, we shall give two examples to illustrate our results.

\section{Preliminaries}
Throughout this paper,  $k$ is an algebraically closed field, $A$ is a finite dimensional $k$-algebra. Denote by ${\rm mod}A$ the category of finitely generated left $A$-modules, $\mathcal{P}({\rm mod}A)$ the full subcategory of mod$A$ consisting of all projective objects in mod$A$, gl.dim$A$ the global dimension of $A$. $D:={\rm Hom}_{k}(-,k)$ denotes the standard duality functor between mod$A$ and mod$A^{\rm op}$. Given a left $A$-module $M$, add$M$ denotes the full subcategory of mod$A$ consisting of all direct summands of finite direct sums of copies of $M$.\\

\noindent {\bf Dimension of a triangulated category.} Let $\mathcal{T}$ be a triangulated category, and $\mathcal{I}, \mathcal{J}$ be subcategories of  $\mathcal{T}$. Denote by $\langle\mathcal{I}\rangle$ the smallest full subcategory of  $\mathcal{T}$ containing $\mathcal{I}$ and closed under finite direct sums, direct summands and shifts.

Denote by $\mathcal{I}\ast\mathcal{J}$ the full subcategory of $\mathcal{T}$  consisting of objects $N$ such that there exists a distinguished triangle $I\rightarrow N\rightarrow J\rightarrow I[1]$ in $\mathcal{T}$, with $I\in \mathcal{I},J\in\mathcal{J} $. Now let $\langle\mathcal{I}\rangle_{0}=0$, $\langle\mathcal{I}\rangle_{1}=\langle\mathcal{I}\rangle$, and inductively $\langle\mathcal{I}\rangle_{n+1}=\langle\langle\mathcal{I}\rangle_{n}\ast\langle\mathcal{I}\rangle\rangle$.
\begin{definition}\cite{Ro2} Let $\mathcal{T}$ be a triangulated category. The dimension of $\mathcal{T}$ is defined as
$${\rm dim}\mathcal{T}:={\rm min}\{n\in \mathbb{N}~|~\exists M \in ob\mathcal{T}~ \text {such that}~ \mathcal{T}=\langle M\rangle_{n+1} \}.$$
\end{definition}
Denote by $\underline{{\rm mod}}A$ the stable module category of $A$. It is defined as follows: the objects of $\underline{{\rm mod}}A$ are the finitely generated $A$-modules, which we will denote by $\underline{M}$ and the homomorphisms in $\underline{{\rm mod}}A$ are given by ${\rm Hom}_{\underline{{\rm mod}}A}(\underline{X},\underline{Y})={\rm Hom}_{A}(X,Y)/\mathcal{I}(X,Y)$, where $\mathcal{I}(X,Y)$ is the subgroup of ${\rm Hom}_{A}(X,Y)$ consisting of the $A$-homomorphisms from $X$ to $Y$ which factor through an injective $A$-module. Denote by $\pi:{\rm mod}A\rightarrow\underline{{\rm mod}}A$ the canonical functor, and $\pi(X)=\underline{X}, \pi(g)=\underline{g}$, where $X$ is an object in ${\rm mod }A$, and $g\in {\rm Hom}_{A}(X,Y)$. The objects $X, Y$ in ${\rm mod}A$ are called $projectively$ $equivalent$ (resp., $injectively$ $equivalent$) if there exist projective (resp., injective) objects $P,Q$ in ${\rm mod}A$ such that $X\oplus P\cong Y\oplus Q$, and denote it by $X\stackrel{P} \sim Y$~(resp., $X\stackrel{I} \sim Y$).

\begin{proposition}\cite[Proposition 1.2.4]{Chen}\label{c-pro} Let $X, Y$ be two objects in {\rm mod}A, then $\pi(X)\cong\pi(Y)$ if and only if $X\stackrel{I} \sim Y$.
\end{proposition}

 Let $A$ be a finite dimensional self-injective algebra. The cosyzygy functor $\Omega^{-1}:\underline{{\rm mod}}A\rightarrow \underline{{\rm mod}}A$ is an equivalence of categories, and a triangulation of $\underline{{\rm mod}}A$ is given by using this functor as a shift and by letting short exact sequences in $\underline{{\rm mod}}A$ correspond to triangles, and hence $\underline{{\rm mod}}A$ is a triangulated category \cite{Hap}. Let $X$ be in mod$A$, and $I(X)$ the injective hull of $X$.  We need the following Lemma.

\begin{lemma}\label{wzg-lemma} Let $A$ be a finite dimensional self-injective algebra and  $X, M $ two objects in ${\rm mod}A$, $d$ a nonnegative integer, then $\underline{X} \in \langle \underline{M}\rangle_{d+1}$ if and only if there exist $r\in \mathbb{N}_{+}, n_{i}\in \mathbb{Z},$ an injective $A$-module $I^{\prime}$, and exact sequences in ${\rm mod}A: 0\rightarrow X_{1}^{i}\rightarrow X^{i}\oplus I(X^{i}_{1})\rightarrow X_{2}^{i}\rightarrow 0$, $i=1,\cdots,r$, where $\underline{X^{i}_{1}}\in \langle \underline{M}\rangle_{d}$, $ \underline{X^{i}_{1}}\in \langle \underline{M}\rangle$ such that $X$ is a direct summand of $\oplus_{i}\Omega^{-n_{i}}(X^{i})\oplus I^{\prime}$
\end{lemma}
\begin{proof} `Necessity'. Since $\underline{X}\in\langle \underline{M}\rangle_{d+1}=\langle\langle \underline{M}\rangle_{d}\ast\langle \underline{M}\rangle\rangle$, there exist $r\in \mathbb{N}_{+}$, $n_{i}\in \mathbb{Z}$, and distinguished triangles $\underline{X^{i}_{1}}\rightarrow \underline{X^{i}}\rightarrow \underline{X^{i}_{2}}\rightarrow \underline{X^{i}_{1}[1]}$, such that $\underline{X}$ is a direct summand of $\oplus_{i}\underline{\Omega^{-n_{i}}(X^{i})}$, where $\underline{X^{i}_{1}}\in \langle \underline{M}\rangle_{d}$, $\underline{X^{i}_{2}}\in\langle \underline{M}\rangle$ for $i=1,\cdots, r$. Then there exist an injective object $I^{\prime}$ in ${\rm mod}A$ such that $X$ is a direct summand of $\oplus_{i}\Omega^{-n_{i}}(X^{i})\oplus I^{\prime}$. Now the assertion follows from \cite[p.65-66]{Rin} that $0\rightarrow X^{i}_{1}\rightarrow X^{i}\oplus I(X^{i}_{1})\rightarrow X^{i}_{2}\rightarrow 0$ is an exact sequence in ${\rm mod}A$, $i=1,\cdots,r$.

`Sufficiency'. Since $0\rightarrow X^{i}_{1}\rightarrow X^{i}\oplus I(X^{i}_{1})\rightarrow X^{i}_{2}\rightarrow 0$ is exact in ${\rm mod}A$, $i=1,\cdots,r$. We have distinguished triangles $\underline{X^{i}_{1}}\rightarrow \underline{X^{i}}\rightarrow \underline{X^{i}_{2}}\rightarrow \underline{\Omega^{-1}(X^{i}_{1})}$, $i=1,\cdots, r$.
Then we have $X^{i}\in \langle \underline{M}\rangle_{d}\ast\langle \underline{M}\rangle \subseteq \langle\underline{M} \rangle_{d+1}$ by the assumption that $\underline{X^{i}_{1}}\in \langle \underline{M}\rangle_{d}, \underline{X^{i}_{2}}\in \langle \underline{M}\rangle$. Then $\underline{X} \in \langle \underline{M}\rangle_{d+1}$ because $X$ is a direct summand of $\oplus_{i}\Omega^{-n_{i}}(X^{i})\oplus I^{\prime}$. \end{proof}

\noindent {\bf Representation dimension.} Let $A$ be a finite dimensional $k$-algebra, a left $A$-module $M$ is a generator~(resp., a cogenerator) in ${\rm  mod}A$ if $A\in {\rm add}M$~(resp.,$DA\in {\rm add}M$). For $A$, the representation dimension rep.dim$A$\cite{Aus} is defined to be
$${\rm rep.dim}A:={\rm inf}\{{\rm gl.dimEnd}_{A}(M)~|~M~\text{ is a generator-cogenerator in}~{\rm mod}A\}.$$

\begin{proposition}\cite[Proposition 3.9]{Rou}\label{Roquier-prop} Let A be a non-semisimple self-injective algebra. Then,
$${\rm dim}(\underline{{\rm mod}}A)+2\leq {\rm rep.dim}A\leq {ll}(A),$$
where ll($A$) denotes the Loewy length of $A$.
\end{proposition}

\noindent {\bf Smash products.} Assume that in this article $A=\oplus_{g\in G}A_{g}$ is a $G$-graded algebra for a finite group $G$ with the identity $e$. For $a\in A$ and $g\in G$, we write $a_{g}$ for the degree-$g$ component of $a$, and $p_{g}$ for the function $G\rightarrow k$ that sends $h$  to $\delta_{h,g}$.

\begin{definition} The $\mathbf{smash~product}$ of $A$ with $G$ is the $k$-algebra $A\#k[G]^*=\oplus_{g\in G}Ap_{g}$ with multiplication given by $ap_{g}\cdot bp_{h}=ab_{gh^{-1}}p_{h}$, $\forall a, b \in A$.
\end{definition}

\noindent {\bf Oppermann dimension.} Let $R = k[x_1,\dots,x_d]$ be the polynomial ring in $d$ variables with coefficients
in $k$ and ${\rm Max} R$ its maximal spectrum, which is the
set of maximal ideals of $R$ equipped with the Zariski topology.
Denote an element of ${\rm Max}R$
by  $\alpha$, and denote by $R_\alpha$  the corresponding localization of $R$, and
$S_\alpha$ the corresponding simple $R$-module.
Denote by  ${\rm fin} R$ be the category of finite length $R$-modules.

A finitely generated $A\otimes_kR$-module $L$ is called a $d$-$dimensional$ $lattice$ if $L$ is projective as an $R$-module. Then there is an exact functor for the $A\otimes_kR$-lattice $L$,
$$
 L\otimes_R-: {\rm fin} R \longrightarrow {\rm mod}A.
$$
It induces the following function, also denoted by $L\otimes_R-:$
$$
  L\otimes_R-:{\rm Ext}^d_R(M,N) \longrightarrow {\rm Ext}^d_A(L\otimes_RM,L\otimes_RN),
  \quad\text{with}\quad (L\otimes_R-)[\epsilon] = [L\otimes_R\epsilon].
$$

One knows that ${\rm Ext}_R^d(S_\alpha,S_\alpha)$ is generated as a $k$-module
by the equivalence class of a long exact sequence of the form
$$
 \epsilon_\alpha :
 \quad 0 \to S_\alpha \to M_1  \to \cdots \to M_d \to S_\alpha \to 0
$$
with $R_\alpha$-modules $M_i$ $(1 \leq i\leq d)$ which are indecomposable and of length 2~(see the argument in \cite[Section 6]{Rin2}).

\begin{definition}
The $d$-dimensional lattice $L$ is called a {\it $d$-dimensional
Oppermann lattice} for $A$ if $O^{d}_{A}(L):=\{\alpha\in {\rm Max} R~|~[L\otimes_R \epsilon_\alpha] \text{ is a non-zero element of } {\rm Ext}_{A}^d(L\otimes_R S_\alpha,L\otimes_R S_\alpha) \}$ is dense in ${\rm Max} R$. The {\it Oppermann dimension} ${\rm Odim}A$ of $A$ is the supremum of $d$ such that there exists a $d$-dimensional Oppermann lattice
$L$ for $A.$
\end{definition}

The following result is due to Oppermann, and it shows that Odim$A$ is always
finite and that one obtains in this way an interesting lower bound for the representation
dimension:
\begin{theorem}\cite{Opp}\label{Opp} Let $A$ be a finite dimensional $k$-algebra which is not semisimple. Then
$${\rm Odim}A+2\leq {\rm rep.dim}A.$$
\end{theorem}

\begin{proposition}\label{Odim} Let $A$ be a finite dimensional $k$-algebra. Then ${\rm Odim}A=0$ if and only if $A$ is of finite representation type.
\end{proposition}
\begin{proof}`Necessity' Noting that $k$ is an algebraically closed field, by \cite[Example 6.2(b)]{Rin2}, we have Odim$A\geq 1$, if $A$ is a representation infinite algebra .

`Sufficiency' If $A$ is semisimple, then Odim$A$=0 by the fact that ${\rm Ext}^1_A(M,N)=0,$  for any $M$, $N$ in mod$A$. Assume that $A$ is a representation-finite non-semisimple algebra, then by Auslander's Theorem ${\rm rep.dim}A=2,$ and hence Odim$A=0$ by Theorem \ref{Opp}. \end{proof}

\section{Representation dimension}

\noindent {\bf Functors of smash products.} \label{triangulated dimensions} Let $A$ be a $G$-graded $k$-algebra for some finite group $G$ with identity $e$. Then there exists a natural embedding of algebras $i: A\hookrightarrow A\#k[G]^*$, sending $a$ to $a\cdot1=\sum _{h\in G} ap_{h}\in A\#k[G]^* $. Denote by $B=A\#k[G]^*$, there is a pull-up functor $F=B\otimes _{A}-:{\rm mod}A\rightarrow  {\rm mod}B $, which is exact since $B$ is  free as a right $A$-module and a push-down functor $i^*: {\rm mod}B\rightarrow {\rm mod}A$ which is induced by the homomorphism $i: A\hookrightarrow A\#k[G]^*$

\begin{definition} \label{separable graded}
Let $A=\oplus_{g\in G}A_{g}$ be a $G$-graded $k$-algebra for some finite group $G$ with identity $e$. $A$ is called $separably$ $graded$ if there is a family of elements $\{x^{g}~|~ g \in G\}$ in the center $A_{e}$ such that
\begin{enumerate}
  \item $\Sigma_{g}x^{g}=1$,

  \item $rx^{g}=x^{hg}r$ for all $r\in A_{h}$ and all $g\in G$.

\end{enumerate}

\end{definition}

 Recall the definition of twisted bimodule. If $\sigma \in {\rm Aut}(A)$ is a $k$-algebra automorphism, and $_{A}M_{A}$ is an $(A,A)$-bimodule, we write $_{1}M_{\sigma}$ for the twisted bimodule, where the left action of $A$ is the same as on $M$, and the right action of $A$ is defined by $ma=m\sigma(a)$, $m\in M, a\in A$.

Let $B=A\#k[G]^*$. One knows that $\{p_g~|~g\in G\}$ is a set of pairwise orthogonal idempotents in $B$ that sum to 1. There is a free right $G$-action on $\{p_g~|~g\in G\}$, given by $p_g\cdot h=p_{gh}$, and this induces a right action of $G$ on $B$.

\begin{lemma} \label{key lemma}
Let $A=\oplus_{g\in G}A_{g}$ be a $G$-graded $k$-algebra for some finite group $G$.
Denote by $B=A\#k[G]^*$, then the following assertions hold.
\begin{enumerate}
\item If $M$ is a projective~(resp., injective) $A$-module, then so is $B\otimes_{A}M$.
\item For an object $M$ in ${\rm mod}A$, we have $B\otimes_{A}\Omega^{-1}M\cong \Omega^{-1}(B\otimes_{A}M)\oplus I$ for an injective $B$-module $I$.
\item Let $Y$ be in ${\rm mod}B$, then $Y$ is a direct summand of $B\otimes_{A}Y$ as $B$-modules.
\item If $A$ is separably graded, then for any $Y \in {\rm mod}A$  , it is a direct summand of $B\otimes_{A}Y$ as $A$-modules.

\end{enumerate}
\end{lemma}
\begin{proof}(1) See \cite[Section 2, Corollary 1]{JJ}.

(2) It can be seen from (1).

(3) By \cite[Lemma 2.2]{Dug}, $B\otimes_{A}B\cong \oplus_{x\in G}(_{1}B_{x})$, hence $B\otimes_{A}Y\cong B\otimes_{A}(B\otimes_{B}Y)\cong(_{B}B\otimes_{A}B)\otimes_{B}Y\cong\oplus_{x\in G}(_{1}B_{x})\otimes_{B}Y\cong Y\bigoplus\oplus_{x\neq e}$($ _{1}B_{x}\otimes_{B}Y)$

(4) See \cite[Theorem 3.1]{JJ}. \end{proof}

\begin{proof}[Proof of Theorem \ref{main-theorem}]
{\bf Claim 1}. Let $_{A}M$ be an object in ${\rm mod}A$ such that $\underline{{\rm mod}}A=\langle\underline{M}\rangle_{n+1}$, then  $\underline{{\rm mod}}B=\langle\underline{B\otimes_{A}M}\rangle_{n+1}$. We proceed by induction on $n$.

$(i)$ $n=0$. Then  $\underline{{\rm mod}}A=\langle\underline{M}\rangle_{1}=\langle\underline{M}\rangle$, and for each $K\in {\rm mod}B$, we have $\underline{_{A}K}\in \langle\underline{M}\rangle$. There exists $r\in \mathbb{N}_{+}, n_{i}\in \mathbb{Z},i=1,\cdots,r,$ and $\underline{K^{\prime}}\in  \underline{{\rm mod}}A$, such that $\underline{K\oplus K^{\prime}}=\underline{\oplus_{i}\Omega^{n_{i}}(M)}$. By Proposition \ref{c-pro}, we have $K\oplus K^{\prime}\stackrel{I}\sim \oplus_{i}\Omega^{-n_{i}}(M)$. By Lemma \ref{key lemma}, $B\otimes_{A}(K\oplus K^{\prime})\stackrel{I}\sim B\otimes_{A}(\oplus_{i}\Omega^{-n_{i}}(M))$ since $B\otimes_{A}-$ is an exact functor and preserves projectivity and injectivity and $B\otimes_{A}(\oplus_{i}\Omega^{-n_{i}}(M))\cong \oplus_{i}\Omega^{-n_{i}}(B\otimes_{A}M)\oplus I$, where $I$ is an injective $B$-module. This implies $B\otimes(K\oplus K^{\prime})\stackrel{I}\sim \oplus_{i}\Omega^{-n_{i}}(B\otimes_{A}M)$ and then $\underline{B\otimes_{A}K}\in \langle\underline{B\otimes_{A}M}\rangle$, and $\underline{K}\in \langle\underline{B\otimes_{A}M}\rangle$ because $K$ is a direct summand of $B\otimes_{A}M$ by Lemma \ref{key lemma}.

$(ii)$ Assume the assertion holds for $n=k$. Let $n=k+1$, and $K$ be in ${\rm mod}B$ such that $\underline{_{A}K}\in \langle\underline{M}\rangle_{k+1+1}$. Then there exist $r\in \mathbb{N}_{+}, n_{i}\in \mathbb{Z},i=1,\cdots,r$, and an injective $A$-mod $I^{\prime}$, and exact sequences in ${\rm mod}A$ as follows.
$$0\longrightarrow K^{i}_{1}\longrightarrow K^{i}\oplus I(K^{i}_{1})\longrightarrow K^{i}_{2}\longrightarrow 0,$$
where $\underline{K^{i}_{1}}\in \langle\underline{M}\rangle_{k+1}$, $\underline{K^{i}_{2}}\in \langle\underline{M}\rangle$, such that $K$ is a direct summand of $\oplus_{i}\Omega^{-n_{i}}(K^{i})\oplus I^{\prime}$. Then $\underline{K}\in \langle\underline{\oplus_{i}K^{i}}\rangle$, and hence $\underline{ B\otimes_{A}K }\in\langle\underline{B\otimes(\oplus K^{i})}\rangle$. On the other hand, we see that
$$\underline{B\otimes_{A}K^{i}_{1}}\longrightarrow \underline{B\otimes_{A}K^{i}}\longrightarrow \underline{B\otimes_{A}K^{i}_{2}}\longrightarrow\underline{\Omega^{-1}(B\otimes_{A}K^{i}_{1})}$$
is an triangle in $\underline{{\rm mod}}B$ for $i=1,\cdots,r$. However, by assumption we have $\underline{B\otimes_{A}K^{i}_{2}}\in \langle\underline{B\otimes_{A}M}\rangle$, $\underline{B\otimes_{A}K^{i}_{1}}\in \langle\underline{B\otimes_{A}M}\rangle_{k+1}$. It implies $\underline{B\otimes_{A}K^{i}}\in \langle\underline{B\otimes_{A}M}\rangle_{k+2}$, $\underline{B\otimes_{A}K}\in \langle\underline{B\otimes_{A}M}\rangle_{k+2}$, and hence $\underline{K}\in \langle\underline{B\otimes_{A}M}\rangle_{k+2}$.

Furthermore, let $A$ be separably graded, we shall show the following claim. It'll imply that ${\rm dim}(\underline{{\rm mod}}A\#k[G]^*)= {\rm dim}(\underline{{\rm mod}} A).$

{\bf Claim 2.} Let $_{B}X$ be in ${\rm mod} B$ such that $\underline{\rm mod}B=\langle\underline{X}\rangle_{n+1}$, we'll show that $\underline{\rm mod}A=\langle\underline{X}\rangle_{n+1}$.
Let $K$ be in ${\rm mod}A$, $B\otimes_{A}K$ is then in ${\rm mod}B$, and hence $\underline{B\otimes_{A}K}\in \underline{{\rm mod}}B=\langle\underline{X}\rangle_{n+1}$. Then by Lemma \ref{wzg-lemma}, there exist $r\in \mathbb{N}_{+}$, $n_{i}\in \mathbb{Z},i=1,\cdots,r$, and exact sequences in ${\rm mod}B$:
$$0\longrightarrow K^{i}_{1}\longrightarrow K^{i}\oplus I(K^{i}_{1})\longrightarrow K^{i}_{2}\longrightarrow 0,$$
where $\underline{K^{i}_{1}}\in \langle\underline{X}\rangle_{n}$,$K^{i}_{2}\in \langle\underline{X}\rangle$, such that $B\otimes_{A}K$ is a direct summand of $\oplus_{i}\Omega^{-n_{i}}(K^{i})\oplus I^{\prime}$ as $B$-modules. Since $A$ is a subalgebra of $B$, the sequences above are also exact in ${\rm mod}A$, and $B\otimes_{A}K$ is also a direct summand of $\oplus_{i}\Omega^{-n_{i}}(Y^{i})\oplus I^{\prime}$ as $A$-modules. Therefore we have distinguished triangles in $\underline{{\rm mod}}A$ as follows
$$ \underline{K^{i}_{1}}\longrightarrow \underline{K^{i}}\longrightarrow \underline{K^{i}_{2}}\longrightarrow \underline{\oplus_{i}\Omega^{-1}(K^{i}_{1}}),$$ and $\underline{K^{i}}\in \langle\underline{X}\rangle_{n+1}$ because $\underline{K^{i}_{1}}\in \langle\underline{K^{i}_{1}}\rangle\in\langle\underline{X}\rangle$. However, $B\otimes_{A}K$ is a direct summand of $\oplus_{i}\Omega^{-1}(K^{i})$, and then $\underline{B\otimes_{A}K}$ is a direct summand of $\underline{\oplus_{i}\Omega^{-1}(K^{i})}$. Since $A$ is separably graded, $K$ is a direct summand of $B\otimes_{A}K$ as $A$-modules by Lemma \ref{key lemma}, we see $\underline{K}\in \langle\underline{X}\rangle_{n+1}$. \end{proof}

\noindent {\bf Oppermann dimension.} Let $R$ be the polynomial ring in $d$ variables with coefficients in $k$ and  $L$ a $d$-dimensional lattice. Recall the notation $O^{d}_{A}(L)$, that is the set of $\alpha \in {\rm Max}R $ such that $[L\otimes_R\epsilon_\alpha]$ is a non-zero element of ${\rm Ext}^{d}_{A}(L\otimes_RS_\alpha,L\otimes_RS_\alpha).$

Let $L$ be a $d$-dimensional lattice for $A$. Then $B\otimes_AL$ is a finitely generated $B\otimes_kR$-module, and it is a projective $R$-module since $B$ is a finitely generated free $A$-module. Hence $B\otimes_AL$ is a $d$-dimensional lattice for $B.$ Conversely, let $L$ be a $d$-dimensional lattice for $B.$ It is clear $i^*(L)$ is a finitely generated $A\otimes_kR$-module, and projective as an $R$-module. That is $i^*(L)$ is a  $d$-dimensional lattice for $A.$

\begin{lemma}\label{lattice} Let Let $A=\oplus_{g\in G}A_{g}$ be a $G$-graded $k$-algebra for some finite group $G$.
Denote by $B=A\#k[G]^*$. Then the following assertions hold.
\begin{enumerate}
  \item Let $L$ be a $d$-dimensional lattice for $B$. Then $O^d_B(L)\subseteq O^d_A(i^*(L))$
  \item Let $A$ be a separably graded $k$-algebra, and $L$ a $d$-dimensional lattice for $A$. Then $O^d_A(L)\subseteq O^d_B(B\otimes_AL).$
\end{enumerate}
\end{lemma}
\begin{proof} (1) Let $\alpha\in O^d_B(L)$. Then $[L\otimes_R\epsilon_\alpha]$ is a non-zero element in ${\rm Ext}^d_{B}(L\otimes_RS_\alpha,L\otimes_RS_\alpha).$ Then $[i^*(L)\otimes_R\epsilon_\alpha]$ is also a non-zero element in  ${\rm Ext}^d_{A}(i^*(L)\otimes_RS_\alpha,i^*(L)\otimes_RS_\alpha)$ by \cite[Theorem 2.2]{JJ}. Indeed, let
$$\rightarrow P_d\rightarrow P_{d-1}\stackrel{f_{d-1}}\rightarrow\cdots\rightarrow P_0\rightarrow L\otimes_RS_\alpha\rightarrow0$$ be a projective resolution of $L\otimes_RS_\alpha$. Denote by $K_d={\rm ker}f_{d-1}$, and $l: K_d \hookrightarrow P_{d-1}$. Then the sequence above induces a projective resolution of $i^*(L)\otimes_RS_\alpha:$
$$\rightarrow i^*(P_d)\rightarrow i^*(P_{d-1})\rightarrow\cdots\rightarrow i^*(P_{0})\rightarrow i^*(L)\otimes_RS_\alpha\rightarrow0,$$ by the fact that $B$ is a finitely generated free $A$-module and $i^*(L_R\otimes S_\alpha)=i^*(L)\otimes_RS_\alpha$. One knows that $[L\otimes_R\epsilon_\alpha]$ can be presented by a $B$-homomorphism $\beta:  K_d\rightarrow L\otimes_RS_\alpha$, and in this case $[i^*(L)\otimes_R\epsilon_\alpha]$  can be presented by $i^*(\beta): i^*(K_d)\rightarrow i^*(L\otimes_RS_\alpha).$ If $[i^*(L)\otimes_R\epsilon_\alpha]=0$, then there exists a $A$-homomorphism $e: i^*(P_{d-1})\rightarrow i^*(L\otimes_RS_\alpha)$ such that $i^*(\beta)=ei^*(l).$ By \cite[Theorem 2.2]{JJ}, $e$ can be lifted to a $B$-homomorphism $\tilde{e}: P_{d-1}\rightarrow L\otimes_RS_\alpha$ with $\tilde{e}(v)=\Sigma_{g\in G}p_{g}e(p_gv)$ for each $v\in P_{d-1}.$ We shall show $\beta=\tilde{e}l$ in the following, and hence $[L\otimes_R\epsilon_\alpha]=0$, contradiction!

Indeed, let $m\in K_d$, then    $\tilde{e}l(m)=\tilde{e}(m)=\Sigma_{g\in G}p_ge(p_gm)$

\ \ \ \ \ \ \ \ \ \ \ \ \ \ \ \ \ \ \ \ \ \ \ \ \ \ \ \ \ \ \ \  $=\Sigma_{g\in G}p_gei^*(l)(p_gm)$

\ \ \ \ \ \ \ \ \ \ \ \ \ \ \ \ \ \ \ \ \ \ \ \ \ \ \ \ \ \ \ \  $=\Sigma_{g\in G}p_gi^*(\beta)(p_gm)$

\ \ \ \ \ \ \ \ \ \ \ \ \ \ \ \ \ \ \ \ \ \ \ \ \ \ \ \ \ \ \ \ $=\Sigma_{g\in G}p_g\beta(p_gm)$

\ \ \ \ \ \ \ \ \ \ \ \ \ \ \ \ \ \ \ \ \ \ \ \ \ \ \ \ \ \ \ \ $=\Sigma_{g\in G}p_gp_g\beta(m)$

\ \ \ \ \ \ \ \ \ \ \ \ \ \ \ \ \ \ \ \ \ \ \ \ \ \ \ \ \ \ \ \ $=(\Sigma_{g\in G}p_g)\beta(m)$

\ \ \ \ \ \ \ \ \ \ \ \ \ \ \ \ \ \ \ \ \ \ \ \ \ \ \ \ \ \ \ \ $=\beta(m).$

(2) Let $\alpha\in O^d_A(L).$ Then $[B\otimes_AL\otimes_R\epsilon_\alpha]$ is a non-zero element in ${\rm Ext}^d_{B}(B\otimes_AL\otimes_RS_\alpha,B\otimes_AL\otimes_RS_\alpha)$ by \cite[Theorem 3.1]{JJ}.
\end{proof}

\begin{proof}[Proof of Theorem \ref{main-theorem-1}] Let $L$ be a $d$-dimensional Oppermann lattice for $B$. Then $i^*(L)$ is a $d$-dimensional Oppermann lattice for $A$ by Lemma \ref{lattice} (1). Then ${\rm Odim}B\leq {\rm Odim}A.$

Let $A$ be a separably graded $k$-algebra, and $L$ a $d$-dimensional Oppermann lattice for $A.$ Then $B\otimes_AL$ is a $d$-dimensional Oppermann lattice for $B$  by Lemma \ref{lattice} (2). Hence ${\rm Odim}A\leq {\rm Odim}B,$ and hence  ${\rm Odim}A={\rm Odim}B.$ \end{proof}

\section{Examples}

\subsection{Exterior algebras}
Let $\wedge V$ be the exterior algebra of an $n$-dimensional vector space $V$ over a field $k$. Then $\wedge V$ is a graded self-injective algebra of Loewy length $n+1$ and then it has a natural $\mathbb{Z}_{n+1}$-grading. By the similar argument as \cite[Theorem 4.2]{Guo},  we have $\wedge V\#k[\mathbb{Z}_{n+1}]^*$ is a graded self-injective algebra of Loewy length $n+1.$ It's not hard to see that if $n+1$ is invertible in the field $k$, then $\wedge V$ is separably graded. By Theorem \ref{main-theorem} and \cite[Theorem 4.1]{Rou}, ${\rm dim}(\underline{{\rm mod}}\wedge V\#k[\mathbb{Z}_{n+1}]^*)={\rm dim}(\underline{{\rm mod}}\wedge V)=n-1$.
Now by Proposition \ref{Roquier-prop} and the argument above, ${\rm rep.dim}(\wedge V\#k[\mathbb{Z}_{n+1}]^*)={\rm rep.dim}(\wedge V)=n+1$.

\subsection{Finite representation type} Let $A$ be a $G$-graded algebra with $G$ a finite group. By Proposition \ref{Odim} and Theorem \ref{main-theorem-1} we have that $A\#k[G]^*$ is of finite representation type if $A$ is.
 It is Corollary 3 in Section 2 of \cite{JJ}. Furthermore, if $A$ is separably graded, then $A\#k[G]^*$ is of finite representation type only if $A$ is. It is Theorem 3.2  of \cite{JJ}.


\begin{thebibliography}{99}

\bibitem{Aus} M. Auslander, Representation dimension of Artin algebras, Queen Mary College Mathematics Notes, London, 1971.
\bibitem{Ber} P. A. Bergh, Representation dimension and finitely generated cohomology, Adv. Math, 2008, 219(1): 389-
400.

\bibitem{Chen} X. W. Chen, Relative singularity categories and Generalized Serre Duality, University of Science and Technology of China, Hefei, 2007~(in Chinese).
\bibitem{Dug} A. Dugas, Periodic resolutions and self-injective algebras of finite type, J. Pure Appl. Algebra, 2010, 214(6): 990-1000.


\bibitem{EHSST} K. Erdmann, M. Holloway, N. Snashall, $\varnothing$. Solberg, R. Taillefer, Support varieties for selfinjective algebras, K-Theory, 2004, 33(1): 67-87.

\bibitem{Guo} J. Y. Guo,  Coverings and truncations of graded self-injective algebras, J. Algebra, 2012, 355(1): 9-34.

\bibitem{Hap} D. Happel, Triangulated Categories in the Representation Theory of Finite Dimensional Algebras, London
Mathematical Society Lecture Note Series 119, Cambridge: Cambridge University Press, 1988.
\bibitem{JJ} A. Jensen, S. J$\varnothing$ndrup, Smash products, group actions, and group graded rings, Math. Scand, 1991, 68: 161-170.
\bibitem{Opp} S. Oppermann, Lower bounds for Auslander's representation dimension, Duke Math. J, 2009, 148: 211-249.

\bibitem{Rin} C. M. Ringel, Represenatation theory of finite dimensional algebra, Lecture notes in Mathematics, 1986, 116: 7-79.
\bibitem{Rin2} C. M. Ringel, On the representation dimension of artin algebras, 2011, preprint, arxiv: 1107.1861v1.
\bibitem{Rou} R. Rouquier, Representation dimension of exterior algebras, Invent. Math, 2006, 165(2): 357-367.
\bibitem{Ro2} R. Rouquier, Dimensions of triangulated categories, J. K-theory, 2008, 1(2): 193-256.


\bibitem{WZ} Q. S. Wu, C. Zhu, Skew group algebras of Calabi-Yau algebras, J. Algebra, 2011, 340(1): 53-76.
\end{thebibliography}
\end{document}